\documentclass[12pt]{amsart}

\usepackage{amssymb}  
\usepackage[all]{xy}

\DeclareFontEncoding{OT2}{}{} 
\newcommand{\textcyr}[1]{%
 {\fontencoding{OT2}\fontfamily{wncyr}\fontseries{m}\fontshape{n}\selectfont #1}
}
\newcommand{\Sha}{{\mbox{\textcyr{Sh}}\!\!}}

\newcommand{\nichts}{\ensuremath{\left.\right.}}
\newcommand{\sigmat}{{\!\! \nichts^\sigma}}
\newcommand{\ssigma}{{\!\! \nichts^\sigma \!}}
\newcommand{\sisigma}{{\!\! \nichts^{\sigma^{-1}} \!}}

\newcommand{\diam}{{\diamondsuit}}

\newcommand{\inv}{{\operatorname{inv}}}
\newcommand{\mytimes}{*{ \!\!\!\!\!\!\!\!\!\!\! \times \!\!\!\!\!\!\!\!\!\!\!}}
\newcommand{\myequals}{*{ \!\!\!\!\!\!\!\!\!\!\! = \!\!\!\!\!\!\!\!\!\!\!}}

\newcommand{\Q}{{\mathbb Q}}

\newcommand{\Z}{{\mathbb Z}}

\newcommand{\Kbar}{{\overline{K}}}
\newcommand{\Abar}{{\overline{A}}}

\newcommand{\oldm}{t}

\newcommand{\calD}{{\mathcal D}}

\newcommand{\umtp}[1]{#1^\times\!/(#1^\times)^2}

\numberwithin{equation}{section}
\newtheorem{theorem}[equation]{Theorem}
\newtheorem{lemma}[equation]{Lemma}
\newtheorem{corollary}[equation]{Corollary}
\newtheorem{proposition}[equation]{Proposition}

\theoremstyle{definition}

\newtheorem{Remark}[equation]{Remark}

\renewcommand{\baselinestretch}{1.1}   

\begin{document}

\title[The Cassels-Tate pairing]{The yoga of the Cassels-Tate pairing}

\author{Tom Fisher}
\address{DPMMS, Centre for Mathematical Sciences, Wilberforce Road,
Cambridge CB3 0WB, UK}
\email{T.A.Fisher@dpmms.cam.ac.uk}

\author{Edward F. Schaefer}
\address{Department of Mathematics and Computer Science,
Santa Clara University, Santa Clara, CA 95053, USA}
\email{eschaefer@scu.edu}

\author{Michael Stoll}
\address{School of Engineering and Science, Jacobs University Bremen,
         P.O.Box 750561, 28725 Bremen, Germany}
\email{m.stoll@jacobs-university.de}

\thanks{The first and second author would like to thank the
hospitality of Jacobs University Bremen. The second author was 
supported by a National Security Agency, Standard Grant and 
a Fulbright Award.}

\date{October 8, 2007}

\begin{abstract}
Cassels has described a pairing on the 2-Selmer group of an elliptic curve
which shares some properties with the Cassels-Tate pairing. In this article,
we prove that the two pairings are the same.
\end{abstract}

\maketitle


\section{Introduction}
\label{Intro}

In \cite{Ca2}, Cassels defined a pairing on the 2-Selmer group of
an elliptic curve over a number field. It shares some
properties with the extension of the Cassels-Tate
pairing to the 2-Selmer group of an elliptic curve over a number
field. He wrote ``It seems highly probable that the two definitions
are always equivalent, but the present writer is no longer an
adept of the relevant yoga.'' (see \cite[p.\ 115]{Ca2}). 
In this article, we prove that the
two pairings are the same.

The Cassels-Tate pairing is an alternating and bilinear pairing
on the Shafarevich-Tate group of an elliptic curve over a number field.
The fact that it is alternating gives information on the structure
of the Shafarevich-Tate group.
For $n\geq 2$,
its extension  from the $n$-torsion of a Shafarevich-Tate group to
an $n$-Selmer group can be used to determine the
image of the $n^2$-Selmer group in the $n$-Selmer group.
This sometimes enables the determination of which elements of
the $n$-Selmer group come from elements of the Mordell-Weil group
and which come from elements of the Shafarevich-Tate group.
The Cassels-Tate pairing is, unfortunately, quite difficult to evaluate
in practice. The pairing defined by Cassels on the 2-Selmer group of an
elliptic curve, however, is quite straightforward to evaluate. So it is
useful to prove that the two pairings are equal on the 2-Selmer group of
an elliptic curve.

In Section~\ref{defns}, we give the Weil-pairing definition and
a new definition of the
Cassels-Tate pairing extended to the $n$-Selmer group of an elliptic curve,
under a hypothesis that is always satisfied for $n$ a prime.
In Section~\ref{TheCaspair} we present
the definition of the pairing defined by Cassels on 2-Selmer groups
of elliptic curves.  In Section~\ref{maindiag} we present a large diagram
and prove it is commutative.
We also discuss why our methods do not easily generalise to $n$-Selmer groups
for $n > 2$.
We use this diagram to prove our main
theorem in Section~\ref{mainthmsection} that the pairing defined by
Cassels is the same as the Cassels-Tate pairing on the 2-Selmer group of
an elliptic curve over a number field. 

Recently, Swinnerton-Dyer \cite{SD} has generalised Cassels'
pairing on the 2-Selmer group, to a
pairing between the $m$-Selmer group and the 2-Selmer group. 
In parallel with the results described above, 
we show that this pairing is again the Cassels-Tate pairing.


\section{Two definitions of the Cassels-Tate pairing}
\label{defns}

Let $E$ be an elliptic curve defined over $K$, a number field. The
Cassels-Tate pairing is a pairing on $\Sha(K,E)$ taking values in $\Q/\Z$.
We refer to \cite{Ca1} for the original definition. 
In the terminology of \cite{PS} this is the homogeneous space
definition.

Let $m, n \ge 2$ be integers. We are interested in the restriction of 
this pairing to the $n$-torsion $\Sha(K,E)[n]$, or more generally
to $\Sha(K,E)[m] \times \Sha(K,E)[n]$. Let $S^n(K,E)$ denote the 
$n$-Selmer group of $E$ over $K$. The group $\Sha(K,E)[n]$ is 
isomorphic to the quotient of $S^n(K,E)$ by the image of 
$E(K)/nE(K)$ under the coboundary map. We write $\langle~,~\rangle_{\rm CT}$
for the extension of the Cassels-Tate pairing to 
$S^m(K,E) \times S^n(K,E)$. By definition this pairing is trivial on
the images of $E(K)/mE(K)$ and $E(K)/nE(K)$.
 
If $M$ is a ${\rm Gal}(\Kbar/K)$-module, then we  
denote $Z^i({\rm Gal}(\Kbar/K),M)$ and $H^i({\rm Gal}(\Kbar/K),M)$ 
by $Z^i(K,M)$ and $H^i(K,M)$, respectively.

We recall an alternative definition of the Cassels-Tate pairing,
called in \cite{PS} the Weil-pairing definition.
For simplicity we assume that the natural map
\begin{equation}
\label{h2map}
 H^2(K,E[n]) \to \prod_{v} H^2(K_v,E[n]), 
\end{equation}
where $v$ runs over all places of $K$, is injective. 
This is known for $n$ a prime \cite[Lemma 5.1]{Ca1}.
(The injectivity does not hold for $E[n]$ replaced by an arbitrary
finite Galois module. See \cite[III.4.7]{Se3} for a counter-example.)
From Section~\ref{TheCaspair} onwards we restrict to the case $n=2$,
so our hypothesis will be automatically satisfied.

Let $a \in S^m(K,E)$ and $a' \in S^n(K,E)$. We apply Galois cohomology
over~$K$ and its completions~$K_v$ to
\[ \xymatrix{0 \ar[r] & E[n] \ar[r] \ar@{=}[d] & E[mn] \ar[r]^{\cdot n} \ar[d]
                      & E[m] \ar[r] \ar[d] & 0 \\
             0 \ar[r] & E[n] \ar[r] & E \ar[r]^{\cdot n} & \ar[r] E \ar[r] & 0
            }
\]
to obtain a commutative diagram
\[ \xymatrix{H^1(K,E[mn]) \ar[r]^{\cdot n} 
               & H^1(K,E[m]) \ar[r] \ar[d] & H^2(K,E[n]) \ar[d] \\
               & \prod_v H^1(K_v,E) \ar[r] & \prod_v H^2(K_v,E[n]) } 
\]
By the above hypothesis, there exists $c \in H^1(K,E[mn])$ with $n c = a$. 
We represent $c$ by a cocycle $\gamma\in Z^1(K,E[mn])$; then
$\alpha = n \gamma \in Z^1(K,E[m])$ represents~$a$.
For each place $v$ of $K$, the cocycle 
${\rm res}_v(\alpha) = \alpha_v$ in
$Z^1(K_v,E(\Kbar_v))$ is a coboundary. So there exists $\beta_v\in
E(\Kbar_v)$ such that $\alpha_v = d\beta_v$ 
(recall $d\beta_v$ is the cocycle $\sigma\mapsto \ssigma \beta_v - \beta_v$). 
Take $Q_v\in  E(\Kbar_v)$ such that $n Q_v = \beta_v$. Consider 
$dQ_v - \gamma_v \in Z^1(K_v,E[n])$, where $\gamma_v$ is the restriction of
$\gamma$. Let $\cup_e$ be the cup product pairing induced by the
Weil pairing from $H^1(K_v,E[n])\times H^1(K_v,E[n])$ to
$H^2(K_v,\mu_n)$. For $s, s'\in H^1(K_v,E[n])$ define $\langle
s,s'\rangle_{{\rm inv}\circ\cup_e,v}$ to be the composition of
$\cup_e$ with the invariant map. 
We define $\langle a, a'\rangle_1=
\sum_v \langle [dQ_v - \gamma_v],a'\rangle_{{\rm inv}\circ\cup_e,v}$.

\begin{proposition}
\label{pair1isCT}
Let $a \in S^m(K,E)$ and $a' \in S^n(K,E)$. 
We have $\langle a, a'\rangle_1 = \langle a, a'\rangle_{\rm CT}$.
\end{proposition}

\begin{proof}
See \cite[Proof of Lemma 4.1]{Ca1} or \cite[\S2.2]{F}. 
\end{proof}

\begin{Remark}
The general form of the Weil-pairing definition, avoiding the 
hypothesis that~(\ref{h2map}) is injective, is given in \cite[p.\ 97]{Mi}.
This variant is used in \cite{PS} to generalise
Proposition~\ref{pair1isCT} to abelian varieties.
\end{Remark}

Let $C$ and $D$ be torsors ({\em i.e.}, principal homogeneous spaces) 
under $E$. A morphism $\pi : D \to C$ is called an {\em $n$-covering}
if $\pi(P+Q) = nP + \pi(Q)$ for all $P \in E$ and $Q \in D$. 
If $C = E$ is the trivial torsor, this coincides with the usual notion
of $n$-covering of~$E$. 
For $Q_1,Q_2 \in D$ we write $Q_1-Q_2$ for the point on~$E$
determined by the fact~$D$ is a torsor under~$E$.
%
Following \cite[Chapter 6]{Stamminger} we define
the coboundary map $\delta_\pi : C(K) \to H^1(K,E[n])$ 
by sending $P \in C(K)$ to the class of $dQ = ( \sigma \mapsto \ssigma Q - Q)$ 
where $Q \in D(\Kbar)$ with $\pi Q =P$. 

In the case $C = E$, there is a standard bijection between the 
$n$-coverings of~$E$ up to $K$-isomorphism, and the Galois cohomology group 
$H^1(K, E[n])$. It is defined as follows. Let $\psi : D \to E$ 
be an isomorphsim of curves over $\Kbar$ with $[n] \circ \psi = \pi$.
Then $\ssigma \psi \circ \psi^{-1}$ is translation by some 
$\xi_\sigma \in E[n]$ and we identify the $K$-isomorphism class of $D$ 
with the class of $\sigma \mapsto \xi_\sigma$ in~$H^1(K, E[n])$.
If $Q \in D(\Kbar)$ with $\pi(Q) = 0$ then 
we can take $\psi : P \mapsto P - Q$,
in which case $D$ is represented by $-dQ$.
Note also that if $C \to E$ is an $m$-covering of~$E$ and $D \to C$ is an 
$n$-covering  of~$C$, then $D \to E$ is an $mn$-covering of~$E$.
If $D \to E$ corresponds to $b \in H^1(K, E[mn])$, 
then $C \to E$ corresponds to $nb \in H^1(K, E[m])$.
 
We give a new definition of the Cassels-Tate pairing,
again under the hypothesis that~(\ref{h2map}) is injective. 
Let $C$ be an $m$-covering of $E$ over $K$ representing $a$. 
By the hypothesis, $a$ is divisible by $n$ in the Weil-Ch\^atelet group.
So there is an $n$-covering $\pi : D \to C$ defined over~$K$.
Let $v$ be a place of $K$. Since
$a$ is trivial in $H^1(K_v,E(\Kbar_v))$, there is a point $P_v\in C(K_v)$.
We define
$\langle a, a'\rangle_2=\sum_v\langle \delta_\pi(P_v),a' 
\rangle_{{\rm inv}\circ\cup_e,v}$.

\begin{proposition}
\label{pair2is1} 
Let $a \in S^m(K,E)$ and $a' \in S^n(K,E)$. 
We have $\langle a,a'\rangle_2 = \langle a, a'\rangle_1$. 
In particular $\langle a,a'\rangle_2$ does not depend on the choice
of the $P_v$.
\end{proposition}

\begin{proof}
Let $R_C\in C(\Kbar)$ and $R_D\in D(\Kbar)$ such that $R_C$ covers
$0$ on $E$ and $R_D$ covers $R_C$. Since $n(dR_D) = dR_C$ represents $-a$, we
can choose $\gamma\in Z^1(K,E[mn])$, as defined above, to be
$-dR_D$. For each place $v$ of $K$ we are given $P_v \in C(K_v)$.
Let $\beta_v = P_v - R_C$, then 
$d \beta_v = -d R_C$; this represents $a \in H^1(K_v,E[m])$.
Take $Q_v \in E(\Kbar_v)$ with $nQ_v= \beta_v$.
Then $dQ_v - \gamma_v  = d(Q_v + R_D)$ and 
$\pi(Q_v + R_D) = \beta_v +R_C = P_v$. 
Hence $\delta_\pi(P_v)$ is represented by the cocycle 
$dQ_v-\gamma_v$ appearing in the definition of $\langle ~,~\rangle_1$.
\end{proof}


\section{The Cassels pairing}
\label{TheCaspair}

In \cite{Ca2}, Cassels defined a bilinear pairing $\langle~,~\rangle_{\rm Cas}$
on $S^2(K,E)$ taking values in $\mu_2$ 
with the following properties. The element $a\in S^2(K,E)$
is in the image of $S^4(K,E)$ precisely when $\langle a,a'\rangle_{\rm Cas}=
+1$ for all $a'\in S^2(K,E)$. For all $a\in S^2(K,E)$ we have $\langle a,
a\rangle = +1$. These are properties of the Cassels-Tate pairing on
a 2-Selmer group as well.

A mild generalisation of Cassels' construction, 
due to Swinnerton-Dyer \cite{SD}, gives a pairing $S^m(K,E) 
\times S^2(K,E) \to \mu_2$. We work with this generalised form of the
pairing, which we continue to denote $\langle~,~\rangle_{\rm Cas}$.
It reduces to Cassels' definition in the case $m=2$.

We prepare to recall the definition of the pairing. 
The group $S^2(K,E)$ is a
subgroup of $H^1(K,E[2])$. Let $\Abar$ be the finite \'{e}tale
algebra that is the Galois module of maps from $E[2]\setminus 0$
to $\Kbar$. Then $\mu_2(\Abar)$ is the Galois module of maps from
$E[2]\setminus 0$ to $\mu_2$. Let $A$ denote the ${\rm
Gal}(\Kbar/K)$-invariants of $\Abar$. Let $E$ be given by
$y^2=F(x)$ where $F(x)=x^3+a_2x^2+a_4x+a_6$ with $a_i\in K$. 
Then $A\cong K[T]/(F(T))$. 
Let $\theta_1, \theta_2,\theta_3$ be the three roots of $F(x)$ in
$\Kbar$.  We have $A\cong
\prod^{\diam}K(\theta_j)$ where $\prod^{\diam}$ denotes taking the
product over one element from each ${\rm Gal}(\Kbar/K)$-orbit of
the set of $\theta_j$'s. Let $T_j= (\theta_j,0) \in E[2] \setminus 0$ 
and define $w:E[2]\to \mu_2(\Abar)$ by $w(P) = ( T_j \mapsto e_2(P,T_j))$. 
Then $w$ induces an injective homomorphism from
$H^1(K,E[2])$ to $H^1(K,\mu_2(\Abar))$, which we also denote $w$. 
Let $r_j$ be the restriction map 
from $H^1(K,\mu_2(\Abar))$ to $H^1(K(\theta_j),\mu_2)$.
Shapiro's Lemma shows that the map $r=\prod^{\diam} r_j$ is 
an isomorphism of $H^1(K,\mu_2(\Abar))$ with 
$\prod^{\diam}H^1( K(\theta_j),\mu_2)$, which we denote $H^1(A,\mu_2)$. 
For each~$j$, we have a Kummer
isomorphism from $H^1(K(\theta_j),\mu_2)$ to $\umtp{K(\theta_j)}$.
This induces an isomorphism, which we denote $k$, from
$H^1(A,\mu_2)$ to $\umtp{A}$. Note that the image of $H^1(K,E[2])$
in $\umtp{A}$, under $k \circ r \circ w$,  
is equal to the kernel of the norm from $\umtp{A}$ to $\umtp{K}$.

We recall the definition of $\langle~,~\rangle_{\rm Cas}$.
Let $a \in S^m(K,E)$ and  $a'\in S^2(K,E)$. 
Let $M=k\circ r\circ w(a')$ be the element of $\umtp{A}$ representing $a'$. 
The element
$a\in S^m(K,E)$ is represented by an $m$-covering $C$ (which Cassels denotes
$\calD_{\Lambda}$) of $E$. Swinnerton-Dyer \cite{SD} shows that
there are rational functions $f_j$ on $C$, 
defined over $K(\theta_j)$, with the following three properties
\renewcommand{\labelenumi}{(\roman{enumi})}
\begin{enumerate}
\item ${\rm div} (f_j) = 2 \calD_j$ 
where $[\calD_j] \mapsto T_j =(\theta_j,0)$ under the 
isomorphism of ${\rm Pic}^0(C)$ and $E$,  
\item each $K$-isomorphism of $K(\theta_i)$ to
$K(\theta_j)$ sending $\theta_i$ to $\theta_j$ sends $f_i$ to $f_j$, 
\item the product $f_1 f_2 f_3$ is a square in $K(C)$.
\end{enumerate}
\renewcommand{\labelenumi}{(\arabic{enumi})}
He then shows that a 2-covering of $C$ may be defined by setting each
$f_j$ equal to the square of an indeterminate. 
In the case $m=2$, Cassels gives an explicit construction of
the $f_j$ (which he denotes $\frac{L_j}{L_0}$)
and this makes it practical to compute the pairing. 
We write $f$ for the element of $A \otimes_K K(C)$ given by $T_j \mapsto f_j$.

Let $v$ be a prime of $K$. Since $C$ represents an element in
$S^m(K,E)$, there is a point $P_v\in C(K_v)$ (which Cassels calls
${\mathfrak C}_v$). For $\gamma_j$, $\delta_j\in \umtp{K_v(\theta_j)}$ we
let $(\gamma_j, \delta_j)_{K_v(\theta_j)}$ denote the 
quadratic Hilbert norm residue symbol.
Let $\Abar_v = A \otimes_{K} \Kbar_v$ and $A_v$ be
its ${\rm Gal}(\Kbar_v/K_v)$-invariants. Then $A_v\cong
\prod^{\diam} K_v(\theta_j)$, where this $\prod^{\diam}$ is taken
over ${\rm Gal}(\Kbar_v/K_v)$-orbits. 
Let $(\gamma,\delta)_{A_v} 
 = \prod^{\diam}(\gamma_j,\delta_j)_{K(\theta_j)_v}$
where $\gamma, \delta \in \umtp{A_v}$ and $\gamma_j, \delta_j$ are their
images in $K_v(\theta_j)^\times/(K_v(\theta_j)^\times)^2$.
Cassels defines 
$\langle a,a'\rangle_{\rm Cas} = \prod_v  (f(P_v),M)_{A_v}.$

\section{The main diagram}
\label{maindiag}

Now let us introduce Figure~\ref{main} which will enable us to prove that
for $a \in S^m(K,E)$, $a' \in S^2(K,E)$ we have 
$\langle a,a'\rangle_{\rm Cas}= \langle a,a'\rangle_2$.
We can define the maps $w$, $r$ and $k$ locally, in an analogous way, and
it will not change the image of $M$, locally. So we will not change our
notation for these maps. 

\begin{equation}
\label{main}
\xymatrix{ 
H^1(K_v,E[2]) \ar[d]_{w} & \mytimes & H^1(K_v,E[2]) \ar[d]_{w} \ar[r]^-{\cup_e} &  H^2(K_v,\mu_2) \ar@{=}[dr] \save[]-<1cm,0.75cm>*\txt{(1)}\restore \\
H^1(K_v,\mu_2(\Abar_v)) \ar[d]_{r}^{\cong} & \mytimes & H^1(K_v,\mu_2(\Abar_v)) \ar[d]_{r}^{\cong} \ar[r]^{\cup_\oldm} &  H^2(K_v,\mu_2(\Abar_v)) \ar[d]_{r}^{\cong} \ar[r]^-{N_*} &  H^2(K_v,\mu_2) \ar[dd]^{\inv} \\
H^1(A_v,\mu_2) \ar[d]_{k}^{\cong} & \mytimes & H^1(A_v,\mu_2) \ar[d]_{k}^{\cong} \ar[r]^{\cup} \save[]+<2cm,0.9cm>*\txt{(2)}\restore 
&  H^2(A_v,\mu_2) \ar[d]^{\prod^\diam \inv_j} 
\save[]+<2cm,0cm>*\txt{(3)}\restore  \\
\umtp{A_v} & \mytimes & \umtp{A_v} \ar[r]^-{\prod^{\diam} 
(~,~)_{K_v(\theta_j)}}  \save[]+<2cm,1cm>*\txt{(4)}\restore 
& \prod^{\diam} \mu_2 \ar[r]^-{\nu} & \mu_2 }
\end{equation}

We identify $\mu_2 \otimes \mu_2 = \mu_2$ 
via $(-1)^p \otimes (-1)^q = (-1)^{pq}$. Since $\mu_2(\Abar_v)$ 
is the Galois module of maps from $E[2]\setminus 0$ to
$\mu_2$, this identification induces a map $\oldm : \mu_2(\Abar_v) \otimes 
\mu_2(\Abar_v) \to \mu_2(\Abar_v)$. 
Let $\cup_\oldm$ be the cup product map via $\oldm$.
Define $N : \mu_2(\Abar_v) \to \mu_2$ by
$(T\mapsto \gamma(T)) \mapsto \prod_T\gamma(T)$, and let 
$N_\ast$ be the map it induces on $H^2$'s. 
Let $r_j$ be the restriction map from $H^2(K_v,\mu_2(\Abar_v))$ 
to $H^2(K_v(\theta_j),\mu_2)$.
In the same was as for the $H^1$'s, Shapiro's Lemma shows that the map
$r=\prod^{\diam} r_j$ is an isomorphism of $H^2(K_v,\mu_2(\Abar_v))$ with
$\prod^{\diam}H^2(K_v(\theta_j),\mu_2)$, which we
denote $H^2(A_v,\mu_2)$.  Let
$\cup_j$ be the cup product map from
$H^1(K_v(\theta_j),\mu_2)\times H^1(K_v(\theta_j),\mu_2)$ to
$H^2(K_v(\theta_j),\mu_2)$ and $\cup = \prod^{\diam}\cup_j$.
Let $\inv : H^2(K_v,\mu_2) \to \mu_2$ be the composition of the invariant
map with the isomorphism of $\frac{1}{2}\Z/\Z$ and $\mu_2$,
and likewise for $\inv_j : H^2(K_v(\theta_j),\mu_2) \to \mu_2$. 
Finally let $\nu : \prod^{\diam}\mu_2 \to \mu_2$ be 
the usual product in~$\mu_2$.
 
\begin{theorem}
\label{diagthm}
The diagram in Figure~\ref{main} is commutative. 
\end{theorem}

We prove this theorem using the following lemmas.

\begin{lemma}
\label{Weilproduct}
Identify $\mu_2\otimes\mu_2 = \mu_2$ as above. 
Then for all $P, Q\in E[2]$ we have
\[
e_2(P,Q)=\prod_{T\in E[2]\setminus 0} e_2(P,T)\otimes e_2(Q,T).\]
\end{lemma}

\begin{proof}
A trivial verification.
\end{proof}

\begin{lemma}
\label{diag1lem}
Diagram (1) in Figure~\ref{main} is commutative.
\end{lemma}

\begin{proof}
Let $\xi, \psi \in H^1(K_v,E[2])$ be represented by cocycles which,
for ease of notation, we also write as $\xi$ and $\psi$.
We have $\xi\cup_e\psi : (\sigma,\tau)
\mapsto e_2(\xi_\sigma,\ssigma\psi_\tau)\in H^2(K_v,\mu_2)$.

Now $w(\xi): \sigma\mapsto (T\mapsto e_2(\xi_\sigma,T))$ for $T\in
E[2] \setminus 0$ and similarly for $w(\psi)$. Thus
\begin{align*}
N_\ast\big( w(\xi) \cup_\oldm w(\psi)\big) : (\sigma,\tau)& \mapsto
N_\ast \big( \oldm \big( (S\mapsto e_2(\xi_\sigma,S))\otimes\ssigma (T\mapsto
e_2(\psi_\tau,T))\big)\big) \\
& =N_\ast \big( \oldm \big( (S\mapsto e_2(\xi_\sigma,S))\otimes (T\mapsto
\ssigma e_2(\psi_\tau,\sisigma T))\big) \big) \\
& =N_\ast \big( \oldm \big( (S\mapsto e_2(\xi_\sigma,S))\otimes (T\mapsto
e_2(\ssigma\psi_\tau,T))\big) \big) \\
& =N_\ast \big( T\mapsto e_2(\xi_\sigma,T)\otimes e_2(\ssigma\psi_\tau,T)\big)
\\
& =\prod_{T\in E[2]\setminus 0} e_2(\xi_\sigma,T)\otimes
e_2(\ssigma\psi_\tau,T) \in \mu_2\otimes\mu_2 . \end{align*}

By Lemma~\ref{Weilproduct} this is the same as $\xi\cup_e\psi$.
\end{proof}

\begin{lemma}
\label{diag2lem}
Diagram (2) in Figure~\ref{main} is commutative
\end{lemma}

\begin{proof}
Let $\xi,\psi \in H^1(K_v,\mu_2(\Abar_v))$. As in the proof of the 
previous lemma, we use the same symbols for 
cocycles representing these classes.
Let $T_j = (\theta_j,0) \in E[2]\setminus 0$.
We must show that $r_j(\xi\cup_\oldm\psi)$ and $r_j(\xi)\cup_j r_j(\psi)$
are equal in $H^2(K_v(\theta_j),\mu_2 \otimes \mu_2)$. We find that
they are represented by cocycles
$(\sigma,\tau) \mapsto\xi_\sigma (T_j)\otimes \ssigma (\psi_\tau) (T_j)$
and
$(\sigma,\tau) \mapsto\xi_\sigma (T_j)\otimes \ssigma (\psi_\tau (T_j))$.
Since $\sigma(T_j)= T_j$ for all 
$\sigma \in {\rm Gal}(\Kbar_v/K_v(\theta_j))$, these cocycles are equal.
\end{proof}

\begin{lemma}
\label{diag3lem}
Diagram (3) in Figure~\ref{main} is commutative.
\end{lemma}

\begin{proof}
We have $\Abar_v = \prod^\diam \overline{K_v(\theta_j)}$ where
$\overline{K_v(\theta_j)} := K_v(\theta_j) \otimes_{K_v} \Kbar_v$. 
Let $N_j$ denote the norm induced by taking the product
over each element in the ${\rm Gal}(\Kbar_v/K_v)$-orbit of $\theta_j$. 
Recall that $\nu : \prod^{\diam}\mu_2 \to \mu_2$ is the usual product 
in $\mu_2$, and let $\nu_*$ be the map it induces on $H^2$'s. Then the map
$N_*$ in Figure~\ref{main} factors as the composite of 
$\prod^\diam N_j$ and $\nu_*$. 

We have the following commutative diagram
\[ 
\xymatrix{ H^2(K_v,\mu_2 (\Abar_v)) \ar[d]^{r} & \myequals & \prod^{\diam} 
H^2(K_v,\mu_2( \overline{K_v(\theta_j)} ) )
\ar[d]^{\prod^{\diam} r_j}  \ar[rr]^-{\prod^\diam N_j} & &
\prod^{\diam} H^2(K_v,\mu_2)  \ar[d]^{\prod^{\diam} {\rm inv}} 
\ar[r]^-{\nu_*} & H^2(K_v,\mu_2) \ar[d]^{{\rm inv}}  \\
H^2(A_v,\mu_2) & \myequals & \prod^{\diam} H^2(K_v(\theta_j),\mu_2) 
\ar[rr]^-{\prod^\diam {\rm inv}_j} \save[]+<3cm,0.9cm>*\txt{(5)}\restore
& & \prod^{\diam} \mu_2 \save[]+<2cm,0.9cm>*\txt{(6)}\restore
\ar[r]^-{\nu} & \mu_2 {\rlap .} }
\]
Diagram (5) commutes by the next lemma. 
That Diagram (6) commutes is obvious. This proves the 
commutativity of Diagram (3).
\end{proof}

\begin{lemma}
\label{newlemma}
Let $X_j$ be the ${\rm Gal}(\Kbar_v/K_v)$-orbit of $T_j$.
There is a commutative diagram
\[ \xymatrix{ H^2(K_v,{\rm Map}(X_j,\mu_{2^\infty})) \ar[r]^-{N_j} 
\ar[d]_{r_j}^{\cong} & H^2(K_v, \mu_{2^\infty}) \ar[d]^{\inv} \\
H^2(K_v(\theta_j),\mu_{2^\infty})  \ar[r]^-{\inv_j} & \Q/\Z {\rlap .}} \]
\end{lemma} 

\begin{proof}
Let $\iota : H^2(K_v, \mu_{2^\infty}) 
\to H^2(K_v,{\rm Map}(X_j,\mu_{2^\infty}))$ be induced by the inclusion
of the constant maps. Then $r_j \circ \iota$ is the restriction
map ${\rm Br}(K_v)[2^\infty] \to {\rm Br}(K_v(\theta_j))[2^\infty]$. By
\cite[\S1 Theorem 3]{Se2} it is multiplication by $d_j$ on the invariants,
where $d_j = [K_v(\theta_j): K_v] = \# X_j$, and therefore surjective.
Since $r_j$ is an isomorphism,
it follows that $\iota$ is surjective. Then for 
$\eta \in H^2(K_v, \mu_{2^\infty})$ we compute
\[ (\inv \circ N_j) (\iota(\eta)) = d_j \inv(\eta) = (\inv_j \circ r_j) 
(\iota (\eta)). \] 
(Alternatively, the definitions in \cite[Chapter III,\S9]{Brown} show that
$N_j \circ r_j^{-1}$ is corestriction, and the lemma then reduces to
a well known property of the invariant maps.)
\end{proof}

\begin{lemma}
\label{diag4lem}
Diagram (4) in Figure~\ref{main} is commutative.
\end{lemma}

\begin{proof}
This is \cite[XIV.2 Prop. 5]{Se1} applied to 
each constituent field of $A_v$.
\end{proof}

Lemmas \ref{diag1lem}, \ref{diag2lem}, \ref{diag3lem}, \ref{diag4lem}
together prove Theorem \ref{diagthm}.
Composing the maps in the last row of~(\ref{main}) 
gives the pairing $(~,~)_{A_v}$ defined at the end of 
Section~\ref{TheCaspair}.
Identifying $\frac{1}{2}\Z/\Z$ with $\mu_2$ we obtain

\begin{corollary}
\label{diagcor}
Let $s,s' \in H^1(K_v,E[2])$. We have
$ \langle s , s' \rangle_{{\rm inv}\circ\cup_e,v}
 = (k \circ r \circ w(s),k \circ r \circ w(s'))_{A_v}$.  
\end{corollary}

\begin{Remark} It would be useful to have an analogue of 
Corollary~\ref{diagcor} for elements of $H^1(K_v,E[n])$
for general $n$ (or at least for $n$ prime).
Lemma~\ref{diag1lem} depends on the equality in 
Lemma~\ref{Weilproduct}, which in turn only works 
for $n=2$. This prevents any obvious generalisation to other
values of $n$. Another difficulty is that 
we use $\mu_2\subset K_v$ in our proofs.
\end{Remark}


\section{The main theorem}
\label{mainthmsection}

Let $C$ be a torsor under $E$, and $f \in A \otimes_K K(C)$ as described
in Section~\ref{TheCaspair}. Let $\pi : D \to C$ be the 2-covering
obtained by setting each $f_j$ equal to the square of an indeterminate.
The following lemma is a variant of \cite[Theorem 2.3]{S}.

\begin{lemma}
\label{delta=f}
We have 
$  (k \circ r \circ w)(\delta_\pi (P)) = f(P) \mod{(A^\times)^2} $
for all $P \in C(K)$, away from the zeroes and poles of the $f_j$.
\end{lemma}

\begin{proof}
Let $Q \in D(\Kbar)$ with $\pi(Q)=P$. It suffices to 
show that $f_j(P ) = k_jr_j w (dQ) \mod{(K(\theta_j)^\times)^2}$. 

We have $r_j w(dQ) = (\sigma \mapsto e_2(\sigmat Q - Q,T_j))$
in $H^1(K(\theta_j),\mu_2)$. The construction of $D$ gives that
$f_j \circ \pi = t_j^2$ for some rational function $t_j$
on $D$, defined over $K(\theta_j)$. 
We claim that $e_2(S,T_j) = t_j(S+X) / t_j(X)$ 
for any  $X \in D(\Kbar)$ for which the numerator
and denominator are well-defined and non-zero. 
Indeed, since the Weil pairing is a geometric construction
we may identify $D$ and $E$ over $\Kbar$. This is an 
identification as torsors under $E$, so the action of $E$ on $D$
is identified with the group law on $E$. Then $\pi$ 
is the multiplication-by-2 map on $E$, and our claim
reduces to the definition of the Weil pairing
in \cite[Chapter III, \S8]{Silverman}.

Putting $S = \sigmat Q - Q$ and $X = Q$ gives
$e_2(\sigmat Q - Q,T_j) = t_j(\sigmat Q)/t_j(Q) = \ssigma
(t_j(Q))/t_j(Q)$ for any $\sigma \in {\rm Gal}(K(\theta_j)/K)$. 
Then $r_j w(dQ) = (\sigma \mapsto \ssigma (t_j(Q))/t_j(Q))$
and $k_j r_j w(dQ) = t_j^2(Q) = f_j \pi (Q) = f_j(P)$.
\end{proof}

As usual we identify $\frac{1}{2}\Z/\Z$ with $\mu_2$.

\begin{theorem}
\label{mainthm}
Let $K$ be a number field and $E$ an elliptic curve over $K$.
Let $a \in S^m(K,E)$ and $a' \in S^2(K,E)$. 
We have $\langle a,a'\rangle_{\rm Cas}= \langle a,a'\rangle_2= \langle a,a'\rangle_1= \langle a,a'\rangle_{\rm CT}$.
\end{theorem}

\begin{proof}
The identification $\langle a,a'\rangle_{\rm Cas}= \langle a,a'\rangle_2$
is immediate from Corollary~\ref{diagcor} and the local analogue of
Lemma~\ref{delta=f}. The other identifications were established in 
Section~\ref{defns}.
\end{proof}



\begin{thebibliography}{19}

\frenchspacing
\renewcommand{\baselinestretch}{1}

\bibitem{Brown}
K.S. \ Brown, {\em Cohomology of groups}, GTM 87, 
Springer-Verlag, New York, 1994. 

\bibitem{Ca1}
 J.W.S.\ Cassels, Arithmetic on curves of genus 1, IV. Proof of the
Hauptvermutung, {\em J.\ reine angew.\ Math.\ } {\bf 211} (1962), 95--112.

\bibitem{Ca2}
  J.W.S.\ Cassels, Second descents for elliptic curves,
 {\em J.\ reine angew.\ Math. } {\bf 494} (1998), 101--127.

\bibitem{F}
T.A. \ Fisher,  The Cassels-Tate pairing and the Platonic solids, 
{\em J. Number Theory  } {\bf{98}} (2003), 105-155.

\bibitem{Mi}
 J.S.\ Milne, {\it Arithmetic duality theorems}, 
Academic Press, Inc., Boston, Mass, 1986.

\bibitem{PS}
B.\ Poonen and M.\ Stoll, The Cassels-Tate pairing on polarized
abelian varieties, {\em Ann.\ of Math.} {\bf 150} (1999), 1109--1149.

\bibitem{S}
E.F. \ Schaefer, 
Computing a Selmer group of a Jacobian using functions on the curve.
{\em Math. Ann.} {\bf 310} (1998), no. 3, 447--471.

\bibitem{Se1}
 J-P.\ Serre, {\it Local fields}, Springer-Verlag, New York, 1979.

\bibitem{Se2}
 J-P.\ Serre, Local class field theory, in J.W.S.\ Cassels and
A.\ Fr\"{o}hlich (eds), {\em Algebraic Number Theory}, 
Academic Press, London, 1967, 129--161.

\bibitem{Se3}
J-P. \ Serre, {\em Galois cohomology}, Springer-Verlag, Berlin, 2002.

\bibitem{Silverman}
J.H. \ Silverman,
{\em The arithmetic of elliptic curves},
GTM 106, Springer-Verlag, New York, 1992. 

\bibitem{Stamminger}
S. \ Stamminger,
{\em Explicit 8-descent on elliptic curves},
PhD thesis, International University Bremen, 2005.

\bibitem{SD}
H.P.F. \ Swinnerton-Dyer, {\em $2^n$-descent on elliptic curves for all $n$},
draught, January 2007.

\end{thebibliography}
\end{document}